\documentclass[reqno,12pt]{amsart}
\usepackage{mathrsfs}
\usepackage{amssymb,latexsym}
\usepackage{amsmath,amsfonts,amsthm}
\usepackage{graphicx}
\usepackage{txfonts}
\usepackage{indentfirst}
\usepackage{fancyhdr}
\usepackage[margin=2.5cm]{geometry}
\allowdisplaybreaks[4]
\usepackage[T1]{fontenc}

\fancypagestyle{plain}{%
 \fancyhead[LO,RE]{\usebox{\mygraphic}}
 
 }
\begin{document}
\setlength{\abovedisplayshortskip}{3mm}
\setlength{\belowdisplayshortskip}{3mm}
\setlength{\abovedisplayskip}{3mm}
\setlength{\belowdisplayskip}{3mm}

\newtheorem{theorem}{Theorem}[section]
\newtheorem{corollary}[theorem]{Corollary}  
\newtheorem{lemma}[theorem]{Lemma}
\newtheorem{definition}[theorem]{Definition}
\newtheorem{proposition}[theorem]{Proposition}
\newtheorem{remark}[theorem]{Remark}
\newtheorem{example}[theorem]{Example}

\newcommand{\bth}{\begin{theorem}}
\newcommand{\eeth}{\end{theorem}}
\newcommand{\ble}{\begin{lemma}}
\newcommand{\ele}{\end{lemma}}
\newcommand{\bco}{\begin{corollary}}
\newcommand{\eco}{\end{corollary}}
\newcommand{\bde}{\begin{definition}}
\newcommand{\ede}{\end{definition}}
\newcommand{\bpr}{\begin{proposition}}
\newcommand{\epr}{\end{proposition}}
\newcommand{\bre}{\begin{remark}}
\newcommand{\ere}{\end{remark}}
\newcommand{\beg}{\begin{example}}
\newcommand{\eeg}{\end{example}}

\newcommand{\beq}{\begin{equation}}
\newcommand{\eeq}{\end{equation}}
\newcommand{\ben}{\begin{equation*}}
\newcommand{\een}{\end{equation*}}
\newcommand{\beqn}{\begin{eqnarray}}
\newcommand{\eeqn}{\end{eqnarray}}
\newcommand{\be}{\begin{eqnarray*}}
\newcommand{\ee}{\end{eqnarray*}}
\newcommand{\ban}{\begin{align*}}
\newcommand{\ean}{\end{align*}}
\newcommand{\bal}{\begin{align}}
\newcommand{\eal}{\end{align}}
\newcommand{\bad}{\aligned}
\newcommand{\ead}{\endaligned}
\newcommand{\lan}{\langle}
\newcommand{\ran}{\rangle}

\newcommand{\na}{\nabla}
\newcommand{\vp}{\varphi}
\newcommand{\La}{\Lambda}
\newcommand{\la}{\lambda}
\newcommand{\Om}{\Omega}
\newcommand{\ta}{\theta}
\newcommand{\fr}{\frac}
\newcommand{\iy}{\infty}
\newcommand{\ve}{\varepsilon}
\newcommand{\pa}{\partial}
\newcommand{\al}{\alpha}
\newcommand{\mr}{\mathbb{R}^n}
\newcommand{\bu}{\bullet}
\newcommand{\si}{\sigma}

\newenvironment{sequation}{\begin{equation}\small}{\end{equation}}
\newenvironment{tequation}{\begin{equation}\tiny}{\end{equation}}

\title{The concavity of $p$-R\'enyi entropy power for doubly nonlinear diffusion equations and $L^p$-Gagliardo-Nirenberg-Sobolev inequalities}
\author
{Yu-Zhao Wang}
\author
{Yan-Mei Wang}
\address{School of Mathematical Sciences, Shanxi University, Taiyuan, 030006, Shanxi, China}
\email{wangyuzhao@sxu.edu.cn}
\email{wangyanmei.mail@qq.com}
\maketitle

\numberwithin{equation}{section}
\numberwithin{theorem}{section}
\setcounter{tocdepth}{2}
\setcounter{secnumdepth}{2}

\begin{abstract}
We prove  the concavity of $p$-R\'enyi entropy power for positive solutions to the doubly nonlinear diffusion equations on $\mathbb{R}^n$ or compact Riemannian manifolds with nonnegative Ricci curvature. As applications, we give  new proofs of the sharp $L^p$-Sobolev inequality and $L^p$-Gagliardo-Nirenberg inequalities on $\mathbb{R}^n$. Moreover, two improvement of $L^p$-Gagliardo-Nirenberg inequalities are derived.
\vspace{2mm}


\vspace{2mm}

\end{abstract}

\section{Introduction and main results}
Entropy power inequality(EPI) states that for any two independent random vectors $X,Y$ on $\mathbb{R}^n$, the entropy power of $X+Y$ is at least the sum of  their entropy powers, that is
\beq\label{EPI}
N(X+Y)\ge N(X)+N(Y),
\eeq
where the entropy power is defined by
\beq\label{EP}
N(X)\doteqdot e^{\frac2nH(X)}
\eeq
 and $"\doteqdot"$ means definition,
\beq\label{Shannonent}
H(X)\doteqdot-\int_{\mathbb{R}^n} u\log u\,dx
\eeq
 is the classic Shannon entropy, $u$ is the probability density of $X$. Equality in \eqref{EPI} holds if and only if the random vectors $X,Y$ are Gaussian with proportional covariances.

This inequality \eqref{EPI} first appeared in Shannon's 1948 historic paper \cite{Shannon} with an incomplete proof, Stam \cite{Stam} gave a complete proof based on Shannon entropy and Fisher information  known as the de Bruijn's identity. Later, many simplified proofs emerged and established connections with other subjects. Lieb \cite{Lieb} proved the EPI via a strengthened Young's inequalty in functional analysis, Dembo-Cover-Thomas \cite{DCT} showed that the similarity with Brunn-Minkowski inequalty in convex geometry, recently, Guo-Shamai-Verd\'u \cite{GSV} obtained a new proof by using of minimum mean-square error(MMSE) in statistic.

On the other hand, in 1985, Costa \cite{Costa} strengthened the EPI for two random vectors in the case where one of these vectors is Gaussian,
\beq\label{Costa1}
N(X+Z_t)\ge (1-t)N(X)+tN(X+Z_1),\quad t\in[0,1]
\eeq
where $Z_t$ is the Gussian with covariance matrix $t{\rm Id}$,
which is equivalent to the concavity of the entropy power of the added Gaussian noise. Namely,
\begin{equation}\label{Concavity}
\frac{d^2}{dt^2}N(X+Z_t)\le0.
\end{equation}
Inequality \eqref{Concavity} is referred to as the \textbf{concavity of entropy
power}.

Later, Demo \cite{Demo} simplified the proof of inequality \eqref{Concavity} based on the Blachman-Stam inequality \cite{Blachman}
\beq
\frac1{I(X+Y)}\le\frac1{I(X)}+\frac1{I(Y)},
\eeq
where $I(X)$ is the Fisher infromation
\beq\label{FI}
I(X)\doteqdot\int_{\mathbb{R}^n}\frac{|\nabla u|^2}{u}\,dx.
\eeq
In \cite{Villani2}, Villani gave a short proof with an exact error relied on the Bakry-\'Emery identities.

The $\gamma$-th R\'enyi entropy of a probability density $u$ in $\mathbb{R}^n$ is
defined by (see, e.g. Cover-Thomas \cite{Cover})
\beq\label{Renyi}
\mathcal{R}_{\gamma}(u)\doteqdot\frac1{1-\gamma}\log\int_{\mathbb{R}^n}u^{\gamma}(x)dx,\quad\gamma\in(0,\infty), \;\gamma\neq1.
\eeq
When $\gamma> 1-\frac2n$, the \textbf{$\gamma$-th R\'enyi entropy power} is given by
\beq\label{RenPow}
\mathcal{N}_{\gamma}(u)\doteqdot\exp\left(\frac{\lambda}{n}\mathcal{R}_{\gamma}(u)\right),
\eeq
where $\lambda=2+n(\gamma-1)>0$.
The R\'enyi entropy for $\gamma= 1$ is defined as the limit of $\mathcal{R}_{\gamma}$
as $\gamma\to1$. It follows directly from definition \eqref{Renyi} that
\beq\label{Shannon}
\mathcal{H}(u)=\lim_{\gamma\to1}\mathcal{R}_{\gamma}(u)=-\int_{\mathbb{R}^n}u(x)\log u(x)dx.
\eeq
Therefore, the Shannon's entropy can be identified with the
R\'enyi entropy of index $\gamma = 1$. In this case, the proposed
R\'enyi entropy power of index $\gamma = 1$, given by \eqref{RenPow}, coincides
with Shannon's entropy power \eqref{EP}.


Recently, G.Savar\'e and G.Toscani \cite{ST} show
that the concavity of entropy power
is a property which is not restricted to Shannon entropy power
 in connection with the heat equation, but it holds for
the $\gamma$-th R\'enyi entropy power \eqref{RenPow}, if we put it in connection
with the solution to the nonlinear diffusion equation
\beq\label{NDE}
\partial_tu=\Delta u^{\gamma}.
\eeq
They show that the $\gamma$-th R\'enyi entropy power defined in \eqref{RenPow} of probability densities in $\mathbb{R}^n$ solving \eqref{NDE} is concave, that is
\beq\label{alphaconva}
\frac{d^2}{dt^2}\mathcal{N}_{\gamma}(u(t))\le0,\quad \gamma > 1-\frac1n.
\eeq
In \cite{WZ} and \cite{WW}, the authors study the concavities of entropy power for positive solutions to the $p$-heat equation and parabolic $p$-Lapaican equation on $\mathbb{R}^n$ and Riemannian manifolds, respectively. In \cite{LL}, Li-Li prove an  entropy power inequality  for positive solutions to the heat equation
of the Witten Laplacian on complete Riemannian manifolds with $CD(K,m)$-condition and on compact manifolds equipped with $(K,m)$-super Ricci flows.

In this paper, we consider the connection between R\'enyi entropy and the doubly nonlinear diffusion equation(DNDE)
\beq\label{DNDE}
\partial_tu=\Delta_p u^{\gamma}
\eeq
on $\mathbb{R}^n$ or compact Riemannian manifold, where $\Delta_p\cdot\doteqdot{\rm div}(|\nabla \cdot|^{p-2}\nabla \cdot)$ is the $p$-Laplacian operator. When $p=2$, DNDE reduces to the equation in \eqref{NDE}.
The  key idea in the study of the nonlinear diffusion equation is to write the equation as a law of mass conservation,
\beq\label{Mass}
\partial_tu+{\rm div}\left(uV\right)=0,
\eeq
which identifies the speed as $V=-\gamma^{p-1}u^{(p-1)(\gamma-1)-1}|\nabla u|^{p-2}\nabla u$ and this in turn allows to write
$V$ as a nonlinear potential flow, $V=-|\nabla v|^{p-2}\nabla v$, this gives for the potential the expression
\beq\label{uv}
 v=\frac{\gamma}{b}u^b,\; b=\gamma-\frac1{p-1}.
\eeq
This potential is just the pressure and the nonlinear speed-pressure relation can also be viewed as the “nonlinear Darcy's law” \cite{WC}.

Thus, the equation \eqref{DNDE} can be equivalently written as
\beq\label{Darcy}
\partial_tu-{\rm div}\left(u|\nabla v|^{p-2}\nabla v\right)=0,
\eeq
and the pressure function $v$ satisfies  the equation
\beq\label{Pressureq}
\partial_tv=bv\Delta_p v+|\nabla v|^p.
\eeq
For a smooth  positive solution $u$ to equation \eqref{DNDE}, define the entropy functional
\beq
\mathcal{E}_b(u)\doteqdot\int_Mu^{b+1}d\mu,\quad  
\eeq
where $M$ is a Riemannian manifold and $d\mu$ is the Riemannian volume measure.
The $p$-R\'enyi entropy is defined by
\beq\label{Renyientropy}
\mathcal{R}_b(u)\doteqdot-\frac1b\log\left(\mathcal{E}_b(u)\right)=-\frac1b\log\left(\int_Mu^{b+1}\,d\mu\right)
\eeq
and the $p$-R\'enyi entropy power
\beq\label{RenyiPower}
\mathcal{N}_{b}(u)\doteqdot\exp\left(\frac{b}{a}\mathcal{R}_{b}(u)\right)
=\Big(\mathcal{E}_b(u)\Big)^{\sigma},\quad
a=-\frac{1}{\sigma}=\frac{nb}{(p-1)nb+p}.
\eeq

Motivated by the works of Villani \cite{Villani2} and Savar\'e-Toscani \cite{ST}, we obtain an explicit second variational formula for $p$-R\'enyi entropy power on Riemannian manifold, when Ricci curvature is nonnegative, $p$-R\'enyi entropy power is concave.

\bth\label{concavity} Let $u$ be a positive solution to \eqref{DNDE} on compact Riemannian manifold $(M,g)$ without boundary, then the $p$-R\'enyi entropy  power  defined in \eqref{RenyiPower} satisfies
\beq\label{concave}
\bad
\frac{d^2}{dt^2}\mathcal{N}_{b}(u)=&pb\sigma\mathcal{E}_b^{\sigma-1}\int_M \left(\Big||\nabla v|^{p-2}\nabla\nabla v-\frac1n(\Delta_pv) a_{ij}\Big|_A^2+|\nabla v|^{2p-4}{\rm Ric}(\nabla v,\nabla v)\right)u^{b+1}\,d\mu\\
&+\sigma(1-\sigma)b^2\mathcal{E}_b^{\sigma-1}\int_M\left|\Delta_pv+\mathcal{I}_b(u)\right|^2u^{b+1}\,d\mu,
\ead\eeq
where $\sigma=-(p-1)-\frac p{nb}$, $a_{ij}$ is the inverse of $A^{ij}=g^{ij}+(p-2)\frac{v^iv^j}{|\nabla v|^2}$, for any $2$-tensor $T$, $|T|^2_A=A^{ij}A^{kl}T_{ik}T_{jl}$ and $\mathcal{I}_b(u)$ is the Fisher information with respect to $p$-R\'enyi entropy $\mathcal{R}_b(u)$
\beq
\mathcal{I}_b(u)\doteqdot\frac{d}{dt}\mathcal{R}_b(u).
\eeq
If Ricci curvature is nonnegative and $b\ge-\frac 1n, b\neq0$, then $t\mapsto \mathcal{N}_b(u)$ is concave.
\eeth
Applying analogous method, the concavity of $p$-R\'enyi entropy power  is valid on $\mathbb{R}^n$.
\bco\label{conRN}
If $u(x,t)$ is a smooth and rapidly decaying nonnegative solution to \eqref{DNDE} on $\mathbb{R}^n$, then the entropy power $\mathcal{N}_b$ is also concave. Moreover,
\beq\label{concave2}
\frac{d^2}{dt^2}\mathcal{N}_{b}(u)=b\sigma\mathcal{E}_b^{\sigma-1}\int_{\mathbb{R}^n} \left[p\Big||\nabla v|^{p-2}\nabla\nabla v-\frac1n(\Delta_pv) a_{ij}\Big|_A^2
+b(1-\sigma)\left(\Delta_pv+\mathcal{I}_b(u)\right)^2\right]u^{b+1}\,dx.
\eeq
In particular, when $p=2$, $b=\gamma-1$ and $\gamma\ge1-\frac1n$, \eqref{concave2} reduces to
\beq\label{concaveform}
\frac{d^2}{dt^2}\mathcal{N}_{\gamma}(u)=2\left(1-\frac2n-\gamma\right)\mathcal{E}_b^{\sigma-1}\int_{\mathbb{R}^n} \left[\Big|\nabla\nabla v-\frac1n(\Delta v) g_{ij}\Big|^2
+\left(\gamma-1+\frac1n\right)\left(\Delta v+\mathcal{I}(u)\right)^2\right]u^{\gamma}\,dx,
\eeq
which give an exact formula in the proof of the concavity of R\'enyi entropy power in \eqref{alphaconva}.
\eco


It is well-know that the fundamental solution to the heat equation on $\mathbb{R}^n$ is
\ben
G(x,t)\doteqdot\frac{1}{(4\pi t)^{\frac n2}}\exp\left\{-\frac{|x|^2}{4t}\right\}.
\een
Since the Shannon's entropy  of $G(x,t)$ equals
\ben
\mathcal{R}(G(x,t))=\frac n2\log(4\pi et).
\een
It follows that the corresponding entropy power is a linear
function of time $t$, i.e.
\ben
\mathcal{N}(G(x,t))=4\pi et,
\een
hence $$\frac{d^2}{dt^2}\mathcal{N}(G(x,t))=0.$$
The concavity property of entropy
power can be rephrased by saying that for all times $t > 0$
the fundamental solution maximizes the second
derivative of the Shannon's entropy power among all possible
solutions to the heat equation.

A natural question is to find the fundamental source-type solution to nonlinear diffusion equation. The first result was found around 1950 by Zel’dovich and Kompaneets and Barenblatt(short BKZ solution) for porous medium equation and fast diffusion equation \eqref{NDE}. For $\gamma>0, \gamma\neq1$, the Barenblatt solution is defined by
\beq\label{Barsol}
U_{\gamma,t}(x)\doteqdot t^{-\alpha}B_{\gamma}(t^{-\frac{\alpha}n}x),\quad \alpha=\frac n{2+(\gamma-1)n}
\eeq
and $B_{\gamma}(\xi)$ is the Barenblatt profile
\ben
B_{\gamma}(\xi)\doteqdot\left\{
                          \begin{array}{ll}
                           (C-\xi^2)_+^{\frac1{\gamma-1}} , & \gamma>1 \\
                            (C+\xi^2)^{\frac1{1-\gamma}}, & 1-\frac 2n<\gamma<1. 
                          \end{array}
                        \right.
\een
In \cite{ST}, Savar\'e-Toscani  showed that the $\gamma$-th R\'enyi entropy power of BKZ solution defined in \eqref{Barsol} is a linear function of time when $\gamma>\frac{ n}{n+2}$
\ben
\mathcal{N}_{\gamma}(U_{\gamma,t})=\mathcal{N}_{\gamma}(U_{\gamma,1})t,
\een
so that
\ben
\frac{d^2}{dt^2}\mathcal{N}_{\gamma}(U_{\gamma,t})=0.
\een
By similar calculations, we can get the source-type solution to the doubly nonlinear diffusion equation in the "good range" (See \cite{VZ})
\ben
\gamma(p-1)>1-\frac pn\Leftrightarrow b>-\frac{p}{n(p-1)}.
\een
The source-type solution is precisely given by
\beq\label{BZK1}
U_{b,t}(x)=t^{-\frac ab}\mathcal{B}_b(t^{-\frac a{nb}}x),
\eeq
where $b=\gamma-\frac1{p-1},\,a=\frac{nb}{nb(p-1)+p}$ and $\mathcal{B}_b(\xi)$ is the Barenblatt profile
\begin{align}\label{Bar}
\mathcal{B}_b(\xi)=\left\{
                   \begin{array}{ll}
                     \Big(C-|\xi|^{\frac p{p-1}}\Big)_+^{\frac1b}, & b>0, \\
                     \Big(C+|\xi|^{\frac p{p-1}}\Big)^{\frac1b}, & -\frac{p}{n(p-1)}<b<0.
                   \end{array}
                 \right.
\end{align}
The parameter $C > 0$ of formula \eqref{Bar} is free and can be uniquely determined
in terms of the unit mass, $\int_{\mathbb{R}^n} U\,dx=1$.

When $b>-\frac{q}{n+q}(q=\frac p{p-1})$, we can also show that the $p$-R\'enyi entropy power
of the source-type solution to the doubly nonlinear diffusion equation is a linear function of $t$,
\beq\label{linear}
\mathcal{N}_b(U_{b,t})=\mathcal{N}_b(U_{b,1})t,
\eeq
that is
\ben
\frac{d^2}{dt^2}\mathcal{N}_{b}(U_{b,t})=0.
\een
The concavity of $p$-R\'enyi entropy in Corollary \ref{conRN} implies that for all $t>0$ the source-type solutions maximize second derivative of the $p$-R\'enyi entropy power among all possible solutions to the doubly nonlinear diffusion equation.

On the other hand, in the same paper \cite{ST}, Savar\'e-Toscani introduced the $\gamma$-weighted Fisher information
\beq
\mathcal{I}_{\gamma}(u)\doteqdot\frac{1}{\int_Mu^{\gamma}\,dx}\int_{\mathbb{R}^n}\frac{|\nabla u^{\gamma}|^2}{u}\,dx,
\eeq
which reduces to the Fisher information \eqref{FI} as $\gamma\to1$. The invariant of scaling and the concavity of R\'enyi entropy power \eqref{alphaconva} imply the following sharp inequality,
\beq\label{REntIso}
\mathcal{N}_{\gamma}(u)\mathcal{I}_{\gamma}(u)\ge \mathcal{N}_{\gamma}(B_{\gamma})\mathcal{I}_{\gamma}(B_{\gamma})=c_{n,\gamma},\quad \gamma>\frac n{n+2}.
\eeq

The sharp inequality \eqref{REntIso} is referred to as the \textbf{isoperimetric inquality for R\'enyi entropy power}. When $\gamma\to1$, \eqref{REntIso} reduces to the linear case
\beq\label{SEntIso}
\mathcal{N}_{1}(u)\mathcal{I}_{1}(u)\ge 2\pi ne=c_{n,1}.
\eeq

Motivated by their works, we introduce the $p$-weighted Fisher information associated with $p$-R\'enyi  entropy
\beq\label{Fisher}
\mathcal{I}_b(u)\doteqdot\frac{b+1}{\gamma\int_Mu^{b+1}\,d\mu}\int_M|\nabla v|^pu\,d\mu=-\frac 1{\int_Mu^{b+1}\,d\mu}\int_M\Delta_pv u^{b+1}\,d\mu.
\eeq
By the identity \eqref{intf1}, we have
\beq\label{DeBruijn}
\frac{d}{dt}\mathcal{R}_b(u)=\mathcal{I}_b(u),\quad t>0.
\eeq
When $p\to2$ and $\gamma\to1$, identity \eqref{DeBruijn} reduces to the classic DeBruijn's identity which connects Shannon's entropy with the Fisher information via heat equation.

A direct computation implies that $\mathcal{I}_b(\mathcal{B}_b)$ is finite if and only if the $q$-th moment of $\mathcal{B}_b$ is finite if $b>-\frac{q}{n+q}$. In this range of exponents, the concavity inequality  \eqref{concave} leads to the following  sharp isoperimetric inequality for $p$-R\'enyi entropy power.

\bth\label{PEntIso}
Let $b>-\frac{q}{n+q}$ and $u$ be the strictly positive and rapidly decaying solution to the doubly nonlinear diffusion equation, then the following sharp inequality holds
\beq\label{idineq}
\mathcal{N}_b(u)\mathcal{I}_b(u)\ge \mathcal{N}_b(\mathcal{B}_b)\mathcal{I}_b(\mathcal{B}_b)=C_{b},
\eeq
where the constant $C_b$ is given by
\beq
C_{1,b}=\left(-\frac{q\gamma}{b}\right)^{p-1}\pi^{\frac p2}n\left(\frac{q(b+1)}{nb+q(b+1)}\right)^{\sigma}\left(\frac{\Gamma(\frac nq+1)\Gamma(-\frac nq-\frac1b)}{\Gamma(\frac n2+1)\Gamma(-\frac1b)}\right)^{\frac pn},\quad -\frac {q}{n+q}<b<0
\eeq
and
\beq
C_{2,b}=\left(\frac{q\gamma}{b}\right)^{p-1}\pi^{\frac p2}n\left(\frac{q(b+1)}{nb+q(b+1)}\right)^{\sigma}\left(\frac{\Gamma(\frac nq+1)\Gamma(\frac1b+1)}{\Gamma(\frac n2+1)\Gamma(\frac nq+\frac1b+1)}\right)^{\frac pn},\quad b>0,
\eeq
where $q=\frac p{p-1}$ and $\sigma=-\left(\frac{p}{nb}+p-1\right)$.

\eeth
Since the inequality \eqref{idineq} is sharp, we can show that an equivalent between inequality \eqref{idineq} and sharp $L^p$-Sobolev inequality when $\gamma=\frac1{p-1}-\frac1n$(i.e. $b=-\frac1n$).
\bth\label{sharpLp}
For any given nonnegative function $w(x)$ such that $w^{p^*}(x)$ is a probability density in $\mathbb{R}^n$, we have
\beq\label{SharpSobolev}
\int_{\mathbb{R}^n}|\nabla w|^p\,dx\ge S_{n,p}
\left(\int_{\mathbb{R}^n}w^{p^*}\,dx\right)^{\frac p{p^*}},\quad p^*=\frac{np}{n-p}>0,
\eeq
where
\ben
S_{n,p}=\left(\frac{p(n-p+1)}{(p-1)(n-p)}\right)^pC_{1,-\frac1n}
=n\pi^{\frac p2}\left(\frac{n-p}{p-1}\right)^{p-1}\left(\frac{\Gamma(\frac nq+1)\Gamma(\frac np)}{\Gamma(\frac n2+1)\Gamma(n)}\right)^{\frac pn},\quad q=\frac p{p-1}
\een
is the sharp $L^p$-Sobolev constant.
\eeth
Moreover, inequality \eqref{idineq} is  equivalent to the
$L^p$-Gagliardo-Nirenberg inequality and we can obtain an improved form of the $L^p$-Gagliardo-Nirenberg inequality.
\bth\label{ImGNI}
{\rm (1)} When $-\frac1n<b<0$, $1<s<\frac n{n-p}$ and $\sigma>1$, we have
\beq\label{PGNI1}
\|w\|_{L^{ps}}\le C_1\|\nabla w\|^{\theta}_{L^p}\|w\|_{L^{(p-1)s+1}}^{1-\theta},
\eeq
where $C_1=[(b+1)\gamma^{p-1}C_{1,b}^{-1} (ps)^p]^{\frac{\theta}p}$ and
\ben\label{theta1}
\theta=\frac1s\frac1{(\sigma+1)\gamma-1}=\frac 1s\frac{n(s-1)}{(p-1)n(1-s)+p((p-1)s+1)}.
\een
Moreover, we get an improvement of the $L^p$-Gagliardo-Nirenberg inequality \eqref{PGNI1},
\begin{align}\label{remainder1}
\|\nabla
w\|_{L^{p}}^{p}\|w\|_{L^{(p-1)s+1}}^{\frac{p(1-\theta)}{\theta}}-C_{1}^{-\frac{p}{\theta}}\|w\|_{L^{ps}}^{\frac{p}{\theta}}
=\left(\frac{1}{ps\gamma}\right)^{p}
\|w\|_{L^{ps}}^{\frac{p}{\theta}}\int_{0}^{\infty}\mathcal{W}(t)dt\ge0,                                                                \end{align}
where
\beq\label{wentropy}
\mathcal{W}_b(t)\doteqdot\frac{\gamma}{b+1}\mathcal{E}_b^{\sigma-1}\int_{\mathbb{R}^n} \left[p\Big||\nabla v|^{p-2}\nabla\nabla v-\frac1n(\Delta_pv) a_{ij}\Big|_A^2
+b(1-\sigma)\left(\Delta_pv+\mathcal{I}_b(u)\right)^2\right]u^{b+1}\,dx.
\eeq
{\rm (2)} When $b>0$, $0<s<1$ and $\sigma<0$, we have
  \beq\label{PGNI2}
\|w\|_{L^{(p-1)s+1}}\le C_2\|\nabla w\|^{\vartheta}_{L^p}\|w\|_{L^{ps}}^{1-\vartheta},
\eeq
where $C_2=[(b+1)\gamma^{p-1}C_{2,b}^{-1} (ps)^p]^{\frac{\vartheta}p}$ and
\ben\label{theta2}
\vartheta=\frac{p}{(1-\sigma)((p-1)s+1)}=\frac 1{(p-1)s+1}\frac{n(1-s)}{n(1-s)+ps}.
\een
Moreover, we get an improvement of the $L^p$-Gagliardo-Nirenberg inequality \eqref{PGNI2},
\begin{align}\label{remainder2}
\|\nabla
w\|_{L^{p}}^{p}\|w\|_{L^{ps}}^{\frac{p(1-\vartheta)}{\vartheta}}-C_{2}^{-\frac{p}{\vartheta}}\|w\|_{L^{(p-1)s+1}}^{\frac{p}{\vartheta}}
=\left(\frac{1}{ps\gamma}\right)^{p}
\|w\|_{L^{(p-1)s+1}}^{\frac{p}{\vartheta}}\int_{0}^{\infty}\mathcal{W}(t)dt\ge0,                                                                \end{align}
where $\mathcal{W}(t)$ is defined in \eqref{wentropy}.
\eeth
\bre
When $p=2$, the results in Theorem \ref{ImGNI} reduce to the improved  $L^2$-Gagliardo-Nirenberg inequality in \cite{DT}.
\ere
\section{The concavity of $p$-R\'enyi entropy power}
Motivated by the methods of Villani \cite{Villani2} and Savar\'e-Toscani \cite{ST},
we need the following identities, which are valid both in $\mathbb{R}^n$ and closed Riemannian manifold with nonnegative Ricci curvature. See the related computations in \cite{WC}.
\begin{lemma}\label{lemma1} For  functions $f$ and $g$, the linearized operator of $p$-Laplacian at point $v$ is defined by
\ben
\mathscr{L}_{p}(f)\doteqdot {\rm div}(|\nabla v|^{p-2}A(\nabla f)),
\een
where $A=g+(p-2)\frac{\nabla v\otimes \nabla v}{|\nabla v|^2}$, then we have
\begin{align}
\label{lpfg}\mathscr{L}_{p}(fg)=&(\mathscr{L}_{p}f)g+f(\mathscr{L}_{p}g)+2|\nabla v|^{p-2}\langle\nabla f,\nabla g\rangle_{A},\\
\label{lpgf}\int_M(\mathscr{L}_{p}f)g\,d\mu=&\int_M f(\mathscr{L}_{p}g)\,d\mu=-\int_M |\nabla v|^{p-2}\langle \nabla f, \nabla g\rangle_{A}\,d\mu,
\end{align}
where $\langle \nabla f, \nabla g\rangle_{A}=A(\nabla f)\cdot\nabla g=\nabla f\cdot A(\nabla g)$.
\end{lemma}
\proof
The proof is a direct result by the definition of $\mathscr{L}_{p}$.
\endproof
\begin{lemma}\label{lemma2}  For $u,v$ in \eqref{uv}, we have
\begin{align}\label{tpv}
\partial_{t}(\Delta_{p}v)=&\mathscr{L}_{p}(\partial_{t}v),\\
\label{tuv}\partial_{t}(uv)=&\frac{b}{p-1}v\mathscr{L}_{p}(uv),
\end{align}
and the $p$-Bochner formula(See in \cite{KoNi} or \cite{WC})
\begin{equation}\label{pbocnner}
\mathscr{L}_{p}|\nabla v|^{p}=p|\nabla v|^{2p-4}(|\nabla\nabla v|^{2}_{A}+{\rm Ric}(\nabla v,\nabla v))+p|\nabla v|^{p-2}\langle\nabla v,\nabla\Delta_{p}v\rangle,
\end{equation}
where
\begin{equation*}
|\nabla\nabla v|_A^2=A^{ij}A^{kl}v_{ik}v_{jl}=\frac{(p-2)^{2}}{4}\frac{|\nabla
v\cdot\nabla|\nabla v|^{2}|^{2}}{|\nabla v|^{4}}+\frac{p-2}{2}\frac{|\nabla|\nabla v|^{2}|^2}{|\nabla v|^{2}}+|\nabla\nabla v|^2.
\end{equation*}
\end{lemma}
\proof
By the definitions of $A$ and the linearized operator $\mathscr{L}_{p}$, we know
\begin{align*}
\partial_{t} (\Delta_{p}v)=&\partial_{t}\left({\rm div}(|\nabla v|^{p-2}\nabla v)\right)={\rm div}\left((p-2)|\nabla v|^{p-4}(\nabla v\cdot\nabla \partial_{t}v)\nabla v+|\nabla v|^{p-2}\nabla\partial_{t} v\right)\nonumber\\
=&{\rm div}\left(|\nabla v|^{p-2}A(\nabla\partial_t v)\right)=\mathscr{L}_{p}(\partial_{t}v).
\end{align*}
For \eqref{tuv},  applying the relationship $u\nabla v=bv\nabla u$ and the definition of $A$, we get \begin{align*}
A(\nabla(uv))
=&\left(1+\frac1b\right)A(u\nabla v)=(p-1)\left(1+\frac1b\right)u\nabla v,
\end{align*}
then
\begin{align*}
\mathscr{L}_{p}(uv)=&{\rm div}\left(|\nabla v|^{p-2}A(\nabla(uv))\right)=(p-1)\left(1+\frac1b\right){\rm div}\left(u|\nabla v|^{p-2}\nabla v\right)=(p-1)\left(1+\frac1b\right)\partial_tu
\end{align*}
and
\begin{align*}
\partial_{t}(uv)=(1+b)v\partial_tu.
\end{align*}
Combining above two equations, we obtain \eqref{tuv}.
\endproof

\bpr\label{intformula}
Let $u$ be a positive solution to \eqref{DNDE} and $v$ satisfies \eqref{Pressureq}, we have
\begin{align}\label{intf1}
\frac{d}{dt}\mathcal{E}_b(u)=&b\int_M(\Delta_{p}v)u^{b+1}\,d\mu=-\frac{b(b+1)}{\gamma}\int_M|\nabla v|^pu\,d\mu=-\frac{b(b+1)}{\gamma}\int_M{u^{-\frac1{p-1}}}{|\nabla u^\gamma|^p}\,d\mu,\\
\label{intf2}
\frac{d^2}{dt^2}\mathcal{E}_b(u)=& \int_M pb\left[|\nabla v|^{2p-4}\Gamma_{2,A}(v)+b(\Delta_{p}v)^2\right]u^{b+1}\,d\mu,
\end{align}
where $\Gamma_{2,A}(v)=|\nabla\nabla v|^{2}_{A}+{\rm Ric}(\nabla v,\nabla v)$.
\epr
\proof

Using the equation \eqref{Darcy} and integrating by parts, we obtain
\ben\bad
\frac{d}{dt}\mathcal{E}_b(u)=&(b+1)\int_Mu^b\partial_tu\,d\mu
=\frac{b(b+1)}{\gamma}\int_Mv{\rm div}\left(u|\nabla v|^{p-2}\nabla v\right)\,d\mu\\
=&-\frac{b(b+1)}{\gamma}\int_M|\nabla v|^pu\,d\mu
=b\int_M(\Delta_{p}v)u^{b+1}\,d\mu,
\ead\een
where we use the identities
\ben
u\nabla v=bv\nabla u,\quad buv=\gamma u^{b+1}.
\een
Applying identities \eqref{intf1}, \eqref{tpv}, \eqref{tuv} and $p$-Bochner formula \eqref{pbocnner}, we have
\ben\bad
\frac{d}{dt}\int_M\Delta_pv(uv)\,d\mu=&\int_M\left[\partial_t(\Delta_pv)(uv)+\Delta_pv\partial_t(uv)\right]\,d\mu\\
=&\int_M\left[\mathscr{L}_{p}(\partial_{t}v)(uv)+\frac b{p-1}\mathscr{L}_{p}(uv)(v\Delta_pv)\right]\,d\mu\\
=&\int_M\left[b\mathscr{L}_p(v\Delta_pv)(uv)+\mathscr{L}_p(|\nabla v|^p)(uv)+\frac b{p-1}\mathscr{L}_{p}(v\Delta_pv)(uv)\right]\,d\mu\\
=&p\int_M\left(|\nabla v|^{2p-4}\Gamma_{2,A}(v)+|\nabla v|^{p-2}\langle\nabla v,\nabla\Delta_{p}v\rangle\right)uvd\mu+\frac{pb}{p-1}\int_M\mathscr{L}_{p}(v\Delta_pv)(uv)\,d\mu.
\ead\een
On the other hand, identities \eqref{lpfg} and \eqref{lpgf} yield
\ben\bad
\int_M\mathscr{L}_{p}(v\Delta_pv)(uv)\,d\mu=&\int_M\left(\mathscr{L}_{p}(v)(\Delta_pv)(uv)+\mathscr{L}_{p}(\Delta_pv)(uv^2)+
2|\nabla v|^{p-2}\langle\nabla \Delta_pv,\nabla v\rangle_A(uv)\right)\,d\mu\\
=&\int_M(p-1)(\Delta_p v)^2(uv)-|\nabla v|^{p-2}\langle\nabla\Delta_pv,\nabla(uv^2)\rangle_A+2|\nabla v|^{p-2}\langle\nabla \Delta_pv,\nabla v\rangle_A(uv)\,d\mu\\
=&\int_M(p-1)(\Delta_p v)^2(uv)-\frac1b|\nabla v|^{p-2}\langle\nabla \Delta_pv,\nabla v\rangle_A(uv)\,d\mu,
\ead\een
where we use the fact
\ben
\nabla(uv^2)=\left(\frac1b+2\right)uv\nabla v.
\een
Combining above formulae, we get
\ben\bad
\frac{d}{dt}\int_Mb(\Delta_pv)u^{b+1}d\mu=&\frac {b^2}{\gamma}\frac{d}{dt}\int_M\Delta_pv(uv)d\mu\\
=&\frac {b^2}{\gamma}\int_M\left(p|\nabla v|^{2p-4}\Gamma_{2,A}(v)+pb(\Delta_pv)^2\right)uvd\mu\\
=&pb\int_M\left(|\nabla v|^{2p-4}\Gamma_{2,A}(v)+b(\Delta_pv)^2\right)u^{b+1}d\mu,
\ead\een
which is \eqref{intf2}.
\endproof

\begin{proof}[\bf{Proof of Theorem \ref{concavity}}]
For a constant $\sigma$, let $\mathcal{N}_{b}(u)=(\mathcal{E}_b(u))^{\sigma}$, by \eqref{intf1} and \eqref{intf2}, we have
\begin{align}\label{2Nbu}
\frac{d^2}{dt^2}\mathcal{N}_{b}(u)
=&\sigma\mathcal{E}_b^{\sigma-2}\left(\mathcal{E}_b\mathcal{E}_b''+(\sigma-1)(\mathcal{E}_b')^2\right)\nonumber\\
=&\sigma\mathcal{E}_b^{\sigma-1}\int_M pb \left(|\nabla v|^{2p-4}\Gamma_{2,A}(v)+b(\Delta_{p}v)^2\right)u^{b+1}\,d\mu
+\sigma(\sigma-1)\mathcal{E}_b^{\sigma-2}\left(\int_M(b\Delta_{p}v)u^{b+1}\,d\mu\right)^2\nonumber\\
=&pb\sigma\mathcal{E}_b^{\sigma-1} \int_M\left[\left(\frac{1}{n}+b\right)(\Delta_pv)^2+\Big||\nabla v|^{p-2}\nabla\nabla v-\frac1n(\Delta_pv) a_{ij}\Big|_A^2+|\nabla v|^{2p-4}{\rm Ric}(\nabla v,\nabla v)\right]u^{b+1}\,d\mu\nonumber\\
&+\sigma(\sigma-1)\mathcal{E}_b^{\sigma-2}\left(\int_M(b\Delta_{p}v)u^{b+1}\,d\mu\right)^2,
\end{align}
where we use the identity
$$|\nabla v|^{2p-4}|\nabla\nabla v|_A^2=\frac1n(\Delta_pv)^2+\Big||\nabla v|^{p-2}\nabla\nabla v-\frac1n(\Delta_pv) a_{ij}\Big|_A^2.$$
When $b>0,\sigma<0$ or $-\frac1n<b<0, \sigma>1$ and $\rm Ric\ge0$, the Cauchy-Schwartz inequality implies that
\begin{align*}
\frac{d^2}{dt^2}\mathcal{N}_{b}(u)\le&\sigma\mathcal{E}_b^{\sigma-1} \int_Mp\left(\frac{1}{nb}+1\right)(b\Delta_pv)^2u^{b+1}\,d\mu+\sigma(\sigma-1)\mathcal{E}_b^{\sigma-2}\left(\int_M(b\Delta_{p}v)u^{b+1}\,d\mu\right)^2\\
\le&\sigma\mathcal{E}_b^{\sigma-1}\left[p \left(\frac{1}{nb}+1\right)+(\sigma-1)\right]\int_M(b\Delta_pv)^2u^{b+1}\,d\mu.
\end{align*}
 Choosing
 \beq\label{sigma}
 \sigma=-\left(\frac{p}{nb}+p-1\right)=-\frac{1}{a},
  \eeq
then \eqref{concave} holds.

In fact, we can obtain an explicit form of $\frac{d^2}{dt^2}\mathcal{N}_{b}(u)$, by the formulae \eqref{2Nbu} and \eqref{sigma}, we have
\begin{align*}
\frac{d^2}{dt^2}\mathcal{N}_{b}(u)=&\sigma(1-\sigma)b^2\mathcal{E}_b^{\sigma-1}\left[\int_M(\Delta v)^2u^{b+1}d\mu-\mathcal{E}_b^{-1}\left(\int_M(\Delta_{p}v)u^{b+1}\,d\mu\right)^2\right]\\
&+pb\sigma\mathcal{E}_b^{\sigma-1} \int_M\left[\Big||\nabla v|^{p-2}\nabla\nabla v-\frac1n(\Delta_pv) a_{ij}\Big|_A^2+|\nabla v|^{2p-4}{\rm Ric}(\nabla v,\nabla v)\right]u^{b+1}\,d\mu\notag\\
=&\sigma(1-\sigma)b^2\mathcal{E}_b^{\sigma-1}\int_M\left|\Delta_pv+\mathcal{I}_b(u)\right|^2u^{b+1}\,d\mu\\
&+pb\sigma\mathcal{E}_b^{\sigma-1}\int_M \left(\Big||\nabla v|^{p-2}\nabla\nabla v-\frac1n(\Delta_pv) a_{ij}\Big|_A^2+|\nabla v|^{2p-4}{\rm Ric}(\nabla v,\nabla v)\right)u^{b+1}\,d\mu,
\end{align*}
where $\mathcal{I}_b(u)=-\frac 1{\int_Mu^{b+1}\,d\mu}\int_M\Delta_pv u^{b+1}\,d\mu$ is the weighted $p$-Fisher information defined in \eqref{Fisher}, when $-\frac1n<b<0$, $\sigma>1$ and $b>0$, $\sigma<0$, thus $\frac{d^2}{dt^2}\mathcal{N}_{b}(u)\le0$.
\end{proof}

\section{R\'enyi entropy and functional inequalities}
\subsection{Isoperimetric inequality for $p$-R\'enyi entropy}

\begin{proof}[\bf Proof of Theorem \ref{PEntIso}]
 By the identity \eqref{DeBruijn}, we get
\ben
\frac{d}{dt}\mathcal{N}_b(u)=\frac ba\mathcal{N}_b(u)\mathcal{I}_b(u).
\een
Set $$\mathcal{Q}_b(u)\doteqdot\mathcal{N}_b(u)\mathcal{I}_b(u),$$
then \eqref{idineq} is equivalent to
\beq
\mathcal{Q}_b(u)\ge\mathcal{Q}_b(\mathcal{B}_b).
\eeq
Since $\frac ba=\frac{nb(p-1)+p}{n}>0$, then by the concavity of $p$-R\'enyi entropy power
\beq\label{decreasing}
\frac{d}{dt}\mathcal{Q}_b(u)=\frac ab\frac{d^2}{dt^2}\mathcal{N}_b(u)\le0,
\eeq
which can be rephrased as the decreasing in time of $\mathcal{Q}_b(u)$.

We can show that $\mathcal{Q}_b$ is invariant with respect to the family of mass-preserving dilations
\ben
D_{\lambda}:u(x)\to D_{\lambda}u(x)\doteqdot\lambda^{-n}u(x/\lambda),\quad\lambda>0.
\een
In fact,
\beq\label{dila}
\mathcal{R}_b(D_{\lambda}u)=\mathcal{R}_b(u)+n\log\lambda,\quad \mathcal{N}_b(D_{\lambda}u)=\lambda^{\frac{bn}{a}} \mathcal{N}_b(u),
\eeq
and
\ben
\mathcal{I}_b(D_{\lambda}u)=\lambda^{-\frac{bn}{a}} \mathcal{I}_b(u),
\een
thus
\beq\label{invariance}
\mathcal{Q}_b(D_{\lambda}u)=\mathcal{Q}_b(u),\quad \lambda>0.
\eeq
An application of  \eqref{dila} to \eqref{BZK1} with $\lambda=t^{\frac a{nb}}$, we can get \eqref{linear}, that is
\ben
\mathcal{N}_b(U_{b,t})=\mathcal{N}_b(U_{b,1})t.
\een
On the other hand, by rescaling the solution
\ben
u_{b,t}(x)=t^{-\frac ab}u(t^{-\frac{an}{b}}x)=D_{t^{\frac{an}{b}}}u,
\een
then
\ben
\mathcal{Q}_b(u_{b,t})=\mathcal{Q}_b(u),
\een
and
\beq\label{longtime}
\lim_{t\to\infty}u_{b,t}(x)=U_{b,1}(x).
\eeq
Combining \eqref{decreasing}, \eqref{longtime} and invariant property \eqref{invariance}, we have
\ben
\mathcal{Q}_b(u)=\mathcal{Q}_b(u_{b,t})\ge \mathcal{Q}_b(U_{b,1})=\mathcal{Q}_b(\mathcal{B}_b).
\een
 For proof of constant $C_b$, we first deduce two  integral formulae by using of the properties of Beta and Gamma functions,
\beq\label{BG1}
\bad
\int_{\mathbb{R}^n}(1-|x|^q)^{\alpha}_+dx=&|\mathbb{S}^{n-1}|\int^1_0\rho^{n-1}(1-\rho^q)^{\alpha}d\rho
=\frac2q\frac{\pi^{\frac{n}{2}}}{\Gamma(\frac n2)}\int^1_0t^{\frac nq-1}(1-t)^{\alpha}dt\\
=&\frac2q\frac{\pi^{\frac{n}{2}}}{\Gamma(\frac n2)}B\left(\frac nq,\alpha+1\right)=\frac2q\frac{\pi^{\frac{n}{2}}}{\Gamma(\frac n2)}\frac{\Gamma(\frac nq)\Gamma(\alpha+1)}{\Gamma(\frac nq+\alpha+1)},
\ead\eeq
and
\beq\label{BG2}
\bad
\int_{\mathbb{R}^n}|x|^q(1-|x|^q)^{\alpha}_+dx=&|\mathbb{S}^{n-1}|\int^1_0\rho^{q+n-1}(1-\rho^q)^{\alpha}d\rho
=\frac2q\frac{\pi^{\frac{n}{2}}}{\Gamma(\frac n2)}\int^1_0t^{\frac nq}(1-t)^{\alpha}dt\\
=&\frac2q\frac{\pi^{\frac{n}{2}}}{\Gamma(\frac n2)}\frac{\Gamma(\frac nq+1)\Gamma(\alpha+1)}{\Gamma(\frac nq+\alpha+2)},\quad q=\frac p{p-1},
\ead\eeq
where $\mathbb{S}^{n-1}$ is the $n-1$ dimensional unit sphere and $|\mathbb{S}^{n-1}|=\frac{2\pi^{\frac n2}}{\Gamma(\frac n2)}$.

(1) The case $b>0$. Set
\ben
D_b\doteqdot\int_{\mathbb{R}^n}(1-|x|^q)^{\frac1b}_+dx=\frac2q\frac{\pi^{\frac{n}{2}}}{\Gamma(\frac n2)}\frac{\Gamma(\frac nq)\Gamma(\frac1b+1)}{\Gamma(\frac nq+\frac1b+1)},
\een
we can obtain $\int_{\mathbb{R}^n}(C-|x|^q)^{\frac1b}_+dx=1$ by choosing
\beq\label{CDB}
C=D_b^{-\frac{bq}{nb+q}}.
\eeq
By \eqref{BG2} and \eqref{CDB}, we have
\beq\label{Mq}
\int_{\mathbb{R}^n}|x|^q\mathcal{B}_b(x)dx=\frac{nb}{nb+q(b+1)}C
\eeq
and
\beq\label{M0}
\int_{\mathbb{R}^n}\mathcal{B}_b(x)^{b+1}dx=\int_{\mathbb{R}^n}(C-|x|^q)\mathcal{B}_b(x)dx=\frac{q(b+1)}{nb+q(b+1)}C.
\eeq
Thus, combining \eqref{Mq} and \eqref{M0} with the definitions of $\mathcal{I}_b$ and $\mathcal{N}_b$, we obtain
\beq\bad
\mathcal{I}_b(\mathcal{B}_b)=&\frac{b+1}{\gamma\int_{\mathbb{R}^n}\mathcal{B}_b(x)^{b+1}dx}
\int_{\mathbb{R}^n}\mathcal{B}_b^{-\frac1{p-1}}|\nabla \mathcal{B}_b^\gamma|^pdx\\
=&\frac{b+1}{\gamma}\left(\frac{q\gamma}{b}\right)^p\frac1{\int_{\mathbb{R}^n}\mathcal{B}_b(x)^{b+1}dx}
\int_{\mathbb{R}^n}|x|^q\mathcal{B}_b(x)dx=\left(\frac{q\gamma}{b}\right)^{p-1}n
\ead\eeq
and
\beq
\mathcal{N}_p(\mathcal{B}_b)=\left(\int_{\mathbb{R}^n}\mathcal{B}_b(x)^{b+1}dx\right)^{\sigma}
=\left(\frac{q(b+1)}{nb+q(b+1)}D_b^{-\frac{bq}{nb+q}}\right)^{\sigma}
=\left(\frac{q(b+1)}{nb+q(b+1)}\right)^{\sigma}D_b^{\frac{p}{n}}.
\eeq
Hence, if $b>0$, the value of the constant $C_{2,b}$ is
\beq
C_{2,b}=\mathcal{N}_p(\mathcal{B}_b)\mathcal{I}_b(\mathcal{B}_b)
=\left(\frac2q\right)^{\frac pn}\left(\frac{q\gamma}{b}\right)^{p-1}\pi^{\frac p2}n\left(\frac{q(b+1)}{nb+q(b+1)}\right)^{\sigma}\left(\frac{\Gamma(\frac nq)\Gamma(\frac1b+1)}{\Gamma(\frac n2)\Gamma(\frac nq+\frac1b+1)}\right)^{\frac pn},
\eeq
where $\sigma=-\left(\frac{p}{nb}+p-1\right)$.

(2) The case $-\frac q{n+q}<b<0$. By analogous computations, we have(or see \cite{VZ})
\beq\label{DB2}
\bad
D_b=&\int_{\mathbb{R}^n}(1+|x|^q)^{\frac1b}dx
=|\mathbb{S}^{n-1}|\int^{\infty}_0\rho^{n-1}(1+\rho^q)^{\frac1b}d\rho\quad\left(s=\frac1{1+\rho^q}\right)\\
=&\frac{|\mathbb{S}^{n-1}|}{q}\int^1_0(1-s)^{\frac nq-1}s^{-\frac1b-\frac nq-1}ds=\frac{|\mathbb{S}^{n-1}|}{q}B\left(-\frac1b-\frac nq,\frac nq\right)=\frac{2\pi^{\frac n2}}{q}\frac{\Gamma(\frac nq)\Gamma(-\frac nq-\frac1b)}{\Gamma(\frac n2)\Gamma(-\frac1{b})},
\ead\eeq
\beq\label{qM}
\bad
\int_{\mathbb{R}^n}|x|^q\mathcal{B}_b(x)dx=&C^{\frac nq+\frac1b+1}\int_{\mathbb{R}^n}|x|^q(1+|x|^q)^{\frac1b}dx=C^{\frac nq+\frac1b+1}\frac{|\mathbb{S}^{n-1}|}{q}\int^1_0(1-s)^{\frac nq}s^{-\frac1b-\frac nq-2}ds\\
=&C^{\frac nq+\frac1b+1}\frac{|\mathbb{S}^{n-1}|}{q}B\left(-\frac1b-\frac nq-1,\frac nq+1\right)=\frac{-nb}{nb+q(b+1)}C,
\ead\eeq
\beq\label{b+1}
\int_{\mathbb{R}^n}\mathcal{B}_b^{b+1}dx=\int_{\mathbb{R}^n}(C+|x|^q)\mathcal{B}_b(x)dx=\frac{q(b+1)}{nb+q(b+1)}C,
\eeq
and
\beq\label{WFIB}
\mathcal{I}_p(\mathcal{B}_b)=\frac{b+1}{\gamma}\left(-\frac{q\gamma}{b}\right)^p\frac1{\int_{\mathbb{R}^n}\mathcal{B}_b(x)^{b+1}dx}
\int_{\mathbb{R}^n}|x|^q\mathcal{B}_b(x)dx=\left(-\frac{q\gamma}{b}\right)^{p-1}n,
\eeq
where $C=D_b^{-\frac{qb}{nb+q}}$.

Combining \eqref{DB2}, \eqref{b+1} and \eqref{WFIB}, we obtain
\beq\label{Cb2}
C_{1,b}=\mathcal{N}_p(\mathcal{B}_b)\mathcal{I}_b(\mathcal{B}_b)
=\left(\frac2q\right)^{\frac pn}\left(-\frac{q\gamma}{b}\right)^{p-1}\pi^{\frac p2}n\left(\frac{q(b+1)}{nb+q(b+1)}\right)^{\sigma}\left(\frac{\Gamma(\frac nq)\Gamma(-\frac nq-\frac1b)}{\Gamma(\frac n2)\Gamma(-\frac1b)}\right)^{\frac pn}.
\eeq
This finished the proof of Theorem \ref{PEntIso}.
\end{proof}
\bre
The formulae \eqref{qM}, \eqref{b+1} and \eqref{WFIB} show that if $b<0$, the $q$-th moment and $p$-weighted Fisher information of the Barenblatt solution $\mathcal{B}_b$ are finite if and only if $b>-\frac q{n+q}$.
\ere
\subsection{Sharp $L^p$-Sobolev inequality}
\begin{proof}[\bf{Proof of Theorem \ref{sharpLp}}]
Let $u(x)$ be a probability density in $\mathbb{R}^n$, if $b>-\frac{q}{n+q}$, inequality \eqref{idineq} can be rewritten as the following form,
\beq\label{spec}
\int_{\{u>0\}}|\nabla v(x)|^pu(x)\,dx=\int_{\{u>0\}}\frac{|\nabla u^{\gamma}(x)|^p}{u^{\frac1{p-1}}(x)}\,dx
\ge \frac{\gamma C_{1,b}}{b+1}\left(\int_{\mathbb{R}^n}u^{b+1}(x)\,dx\right)^{p+\frac{p}{nb}}.
\eeq
Consider the special case $b=-\frac1n$, this leads to
\ben
p+\frac p{nb}=0,\quad {\rm and}\quad\mathcal{N}_{-\frac1n}(u)=\int_{\mathbb{R}^n}u^{1-\frac1n}(x)dx.
\een
Note that the restriction $p<n$ implies $-\frac1n>-\frac q{n+q}$, then for $b=-\frac1n$, we have
\beq\label{Sharp0}
\int_{\mathbb{R}^n}\frac{|\nabla u^{\gamma}(x)|^p}{u^{\frac1{p-1}}(x)}\,dx
\ge \frac{\gamma}{b+1}C_{1,b}.
\eeq
The substitution $u(x)=w^{p^*}(x)$($p^*=\frac{np}{n-p}$) yields
\beq\label{Sharp1}
\int_{\mathbb{R}^n}\frac{|\nabla u^{\gamma}(x)|^p}{u^{\frac1{p-1}}(x)}\,dx
=\left(\frac{p(n-p+1)}{(p-1)(n-p)}\right)^p\int_{\mathbb{R}^n}|\nabla w|^p\,dx.
\eeq
Thus, combining \eqref{Sharp0}, \eqref{Sharp1} and \eqref{Cb2}, we have
 \beq
 \int_{\mathbb{R}^n}|\nabla w|^p\,dx\ge n\pi^{\frac p2}\left(\frac{n-p}{p-1}\right)^{p-1}\left(\frac{\Gamma(\frac np)\Gamma(\frac nq+1)}{\Gamma(\frac n2+1)\Gamma(n)}\right)^{\frac pn}.
 \eeq
The scaling argument shows that if $w^{p^*}(x)$ is a probability density function, then $w$ satisfies the sharp $L^p$-Sobolev inequality \eqref{SharpSobolev}.  In other words, optimal $L^p$ Sobolev inequality is a consequence of the concavity of $p$-R\'enyi entropy power with the parameter $\gamma=\frac{1}{p-1}-\frac1n$.
\end{proof}

\subsection{$L^p$-Gagliardo-Nirenberg inequality}
In \cite{DT}, the authors observe that the isoperimetric inequality
\beq\label{NIu}
\mathcal{N}_{\gamma}(u)\mathcal{I}_{\gamma}(u)\ge \mathcal{N}_{\gamma}(\mathcal{B}_{\gamma})\mathcal{I}_{\gamma}(\mathcal{B}_{\gamma})=c_{\gamma}
\eeq
for the R\'enyi entropy power with nonlinear diffusion equation
\ben
\partial_tu=\Delta u^{\gamma}
\een
is equivalent to one of the following Gagliardo-Nirenberg inequalities(GNI):
\beq\label{GNI1}
\|w\|_{L^{2s}}\le C_1\|\nabla w\|^{\theta}_{L^2}\|w\|_{L^{s+1}}^{1-\theta},\quad 1-\frac1n<\gamma<1
\eeq
and
\beq\label{GNI2}
\|w\|_{L^{s+1}}\le C_2\|\nabla w\|^{\vartheta}_{L^2}\|w\|_{L^{2s}}^{1-\vartheta},\quad \gamma>1,
\eeq
where
\ben
\theta=\frac 1s\frac{n(s-1)}{n(1-s)+2(s+1)},\quad\vartheta=\frac{1}{s+1}\frac{n(1-s)}{n(1-s)+2s}
\een
and
\ben
u=w^{2s}, \quad s=\frac1{2\gamma-1}.
\een
Straightforward computations imply \eqref{NIu} can be brought into the form
\beq
\left(\int_Mu\,dx\right)^{(\sigma+1)\gamma-1}\le \gamma^2c^{-1}_{\gamma} \left(\int_{\mathbb{R}^n}|\nabla u|^2u^{2\gamma-3}\,dx\right)\left(\int_{\mathbb{R}^n}u^{\gamma}\,dx\right)^{\sigma-1},\quad \sigma=\frac{2}{n(1-\gamma)}-1.
\eeq

Motivated by this observation and the isoperimetric inequality \eqref{idineq}, we expect to obtain $L^p$-Gagliardo-Nirenberg inequalities. First, \eqref{idineq} can be rewritten as
\beq\label{PGNS2}
\left(\int_Mudx\right)^{(\sigma+1)\gamma-1}\le (b+1)\gamma^{p-1}C_b^{-1} \left(\int_{\mathbb{R}^n}|\nabla u|^pu^{pb+1-p}\,dx\right)\left(\int_{\mathbb{R}^n}u^{b+1}\,dx\right)^{\sigma-1},
\eeq
where $\sigma=-\left(\frac{p}{nb}+p-1\right)$, $C_b=C_{1,b}$ if $-\frac q{n+q}<b<0$, $C_b=C_{2,b}$ if $b>0$.
Second, set
\beq\label{wps}
u=w^{ps}, \quad s=\frac1{pb+1},
\eeq
and put \eqref{wps} into \eqref{PGNS2},  we have
\beq\label{PGNS3}
\left(\int_Mw^{ps}d\mu\right)^{(\sigma+1)\gamma-1}\le C \left(\int_{\mathbb{R}^n}|\nabla w|^p\,dx\right)\left(\int_{\mathbb{R}^n}w^{(p-1)s+1}\,dx\right)^{\sigma-1},
\eeq
where $C=(b+1)\gamma^{p-1}C_b^{-1} (ps)^p$.
If $-\frac1n<b<0$ and $\sigma>1$, then \eqref{PGNS3} is equivalent to
\beq\label{PGNS1}
\|w\|_{L^{ps}}\le C_1\|\nabla w\|^{\theta}_{L^p}\|w\|_{L^{(p-1)s+1}}^{1-\theta},
\eeq
where $C_1=[(b+1)\gamma^{p-1}C_{1,b}^{-1} (ps)^p]^{\frac{\theta}p}$ and
\beq\label{theta1}
\theta=\frac1s\frac1{(\sigma+1)\gamma-1}=\frac 1s\frac{n(s-1)}{(p-1)n(1-s)+p((p-1)s+1)},\quad s>1.
\eeq
If $b>0$ and $\sigma<0$,  then \eqref{PGNS3} is equivalent to
\beq\label{PGNS4}
\|w\|_{L^{(p-1)s+1}}\le C_2\|\nabla w\|^{\vartheta}_{L^p}\|w\|_{L^{ps}}^{1-\vartheta},
\eeq
where $C_2=[(b+1)\gamma^{p-1}C_{2,b}^{-1} (ps)^p]^{\frac{\vartheta}p}$ and
\beq\label{theta2}
\vartheta=\frac{p}{(1-\sigma)((p-1)s+1)}=\frac 1{(p-1)s+1}\frac{n(1-s)}{n(1-s)+ps},\quad 0<s<1.
\eeq
Furthermore, according to identity \eqref{intf1}, we define a quality
\ben
J_b(u)\doteqdot
-\frac{\gamma}{b(b+1)}\mathcal{E}_b^{\sigma-1}(u)\frac{d}{dt}\mathcal{E}_b(u)
=\frac{\gamma}{b+1}\mathcal{E}_b^{\sigma}(u)\mathcal{I}_b(u),
\een
then identity \eqref{concave2} implies that
\ben
\frac{d}{dt}J_b(u)=\frac{\gamma}{\sigma b(b+1)}\frac{d^2}{dt^2}\mathcal{N}_b(u)=-\mathcal{W}_b(u),
\een
where
\beq\label{went}
\mathcal{W}_b(t)\doteqdot\frac{\gamma}{b+1}\mathcal{E}_b^{\sigma-1}\int_{\mathbb{R}^n} \left[p\Big||\nabla v|^{p-2}\nabla\nabla v-\frac1n(\Delta_pv) a_{ij}\Big|_A^2
+b(1-\sigma)\left(\Delta_pv+\mathcal{I}_b(u)\right)^2\right]u^{b+1}\,dx.
\eeq
Thus, for all $t\ge0$, we get
\beq\label{IGN}
J_p(u_{0})=J_p(u_{\infty})+\int_{0}^{\infty} \mathcal{W}_b(t)dt,
\eeq
where $u_0$ and $u_{\infty}$ are the value of $u(t,x)$ at time $0$ and $\infty$
respectively. \\
(1)When $-\frac1n<b<0$, $\sigma>1$ and $1<s<\frac n{n-p}$, set $u_0(x)=\frac{w^{ps}(x)}{\|w\|^{ps}_{L^{ps}}}$, then we have
\begin{align*}
J_p(u_{0})=&\mathcal{E}_b^{\sigma-1}(u_0)\int_{\mathbb{R}^n}{u_0^{-\frac1{p-1}}}{|\nabla u_0^\gamma|^p}\,dx=\left[\int_{\mathbb{R}^n}\left(\frac{w(x)}{\|w\|_{L^{ps}}}\right)^{ps(b+1)}dx\right]^{\sigma-1}
\int_{\mathbb{R}^n}(ps\gamma)^p\frac{|\nabla w|^p}{\|w\|^p_{L^{ps}}}\,dx\\
=&(ps\gamma)^p\frac{\|\nabla w\|^p_{L^{p}}\|w\|_{L^{ps(b+1)}}^{ps(b+1)(\sigma-1)}}{\|w\|_{L^{ps}}^{ps(b+1)(\sigma-1)+p}}
=(ps\gamma)^p\frac{\|\nabla w\|^p_{L^{p}}\|w\|_{L^{ps(b+1)}}^{p(1-\frac1{\theta})}}{\|w\|_{L^{ps}}^{\frac{p}{\theta}}},
\end{align*}
where
$
\theta=\frac{1}{s(b+1)(\sigma-1)+1}=\frac1s\frac1{(\sigma+1)\gamma-1}
$
is the same with the definition in \eqref{theta1}, set
\ben
C_{1}^{-\frac{p}{\theta}}\doteqdot\left(\frac{1}{ps\gamma}\right)^{p}J_{b}(u_{\infty}),
\een
then we obtain an improvement of $L^p$-Gagliardo-Nirenberg inequality in \eqref{remainder1} by formula \eqref{IGN}.\\
(2)When $b>0$, $\sigma<0$ and $0<s<1$, then
\begin{align*}
J_p(u_{0})=&(ps\gamma)^p\left[\int_{\mathbb{R}^n}\left(\frac{w(x)}{\|w\|_{L^{ps}}}\right)^{ps(b+1)}dx\right]^{\sigma-1}
\int_{\mathbb{R}^n}\frac{|\nabla w|^p}{\|w\|^p_{L^{ps}}}\,dx\\
=&(ps\gamma)^p\left[\frac{\|w\|_{L^{ps}}^{ps(b+1)}}{\|w\|_{L^{ps(b+1)}}^{ps(b+1)}}\right]^{1-\sigma}\|\nabla
w\|_{L^{p}}^{p}\|w\|_{L^{ps}}^{-p}
=\left(ps\gamma\right)^p\frac{\|\nabla
w\|_{L^{p}}^{p}\|w\|_{L^{ps}}^{\frac{p(1-\vartheta)}{\vartheta}}}{\|w\|_{L^{ps(b+1)}}^\frac{p}{\vartheta}},   \end{align*}
where $\vartheta=\frac{1}{s(b+1)(1-\sigma)}=\frac{p}{(1-\sigma)((p-1)s+1)}$, by  \eqref{IGN},  we obtain an improvement of $L^p$-Gagliardo-Nirenberg inequality in \eqref{remainder2}.

\section*{Acknowledgements}
This work is partially supported by the National Natural Science Foundation of China (NSFC Grant No. 11701347).

\end{document}